\documentclass[leqno]{article}

\title{Strong converse inequalities for the weighted multivariate Bernstein-Durrmeyer operator on the simplex via multipliers}
\author{Borislav R. Draganov}
\date{}

\usepackage{amsmath,amsthm,amsfonts}
\usepackage[english]{babel}
\usepackage{fancyhdr}

\allowdisplaybreaks[1]

\newtheorem{thm}{Theorem}[section]
\newtheorem{prop}[thm]{Proposition}

\newtheorem{lem}[thm]{Lemma}

\theoremstyle{remark}
\newtheorem{rem}[thm]{\bf Remark}

\numberwithin{equation}{section}

\newcommand{\thmref}[1]{Theorem~\ref{#1}}
\newcommand{\lemref}[1]{Lemma \ref{#1}}
\newcommand{\propref}[1]{Proposition \ref{#1}}

\newcommand{\N}{\mathbb{N}}
\newcommand{\R}{\mathbb{R}}

\begin{document}

\maketitle
\bigskip

\thispagestyle{fancy}
\fancyhf{}
\renewcommand{\headrulewidth}{0pt}
\lhead{\textit{Ann. Sofia Univ., Fac. Math and Inf.} \textbf{105} (2018), 55--73.}
\renewcommand{\footrulewidth}{0.3pt}
\lfoot{\copyright\,2018. This manuscript version is made available under the CC-BY-NC-ND 4.0 license \texttt{http://creativecommons.org/licenses/by-nc-nd/4.0/}}

\begin{abstract}
\noindent
It is demonstrated that multiplier methods naturally yield better constants in strong converse inequalities for the Bernstein-Durrmeyer operator. The absolute constants obtained in some of the inequalities are independent of the weight and the dimension. The estimates are stated in terms of the $K$-functional that is naturally associated to the operator.
\end{abstract}

\bigskip
\noindent
{\footnotesize \leftskip25pt \rightskip25pt {\sl  AMS} {\it classification}: 41A10, 41A25, 41A27, 41A35, 41A36, 41A63, 42B15.\\[2pt]
\par}
\bigskip

\section{A characterization of the rate of approximation of the Bernstein-Durrmeyer operator}

For $x=(x_1,\dotsc,x_d)\in\R^d$ we set $|x|:=\sum_{i=1}^d |x_i|$. Let $S$ be the standard simplex in $\R^d$ given by
$$
S:=\{(x_1,\dotsc,x_d)\in\R^d : x_i\ge 0,\ i=1,\dotsc,d,\ |x|\le 1\}.
$$
The Jacobi weights on $S$ are defined by
\begin{equation}\label{defw}
w_{\alpha}(x):=x_1^{\alpha_1}\dotsm x_d^{\alpha_d}(1-|x|)^{\alpha_{d+1}},\quad \alpha_i>-1,\ i=1,\dotsc,d+1.
\end{equation}
We have set $\alpha:=(\alpha_1,\dotsc,\alpha_{d+1})$. For $p\in [1,\infty)$ and a Jacobi weight $w_\alpha$ we consider the space $L_{p,w_\alpha}(S)$ of Lebesgue measurable functions $f$ defined on $S$ such that
$$
\|f\|_{p,w_\alpha}:=\left(\int_S |f(x)|^p w_\alpha(x)\,dx\right)^{1/p}<\infty.
$$
Let, as usual, $L_\infty(S)$ denote the space of the essentially bounded Lebesgue measurable functions on $S$, equipped with the sup-norm on $S$. For brevity we set $L_{\infty,w_\alpha}(S):=L_\infty(S)$ and $\|f\|_{\infty,w_\alpha}:=\mathrm{ess\,sup}_{x\in S}|f(x)|$. 

We proceed to the definition of the multivariate Bernstein-Durrmeyer operators with Jacobi weights given by Ditzian \cite{Di:MBD}. For $n\in\N_0$ and $\mathbf{k}=(k_1,\dotsc,k_d)\in\N_0^d$ with $|\mathbf{k}|\le n$ we define the polynomials
$$
p_{n,\mathbf{k}}(x):=\frac{n!}{k_1!\dotsm k_d!(n-|\mathbf{k}|)!}\,\prod_{i=1}^d x_i^{k_i}(1-|x|)^{n-|\mathbf{k}|}. 
$$
The Jacobi-weighted Bernstein-Durrmeyer operators on $L_{p,w_\alpha}(S)$ are defined by
$$
M_{n,\alpha} f(x):=\sum_{|\mathbf{k}|\le n} p_{n,\mathbf{k}}(x)\left(\int_S p_{n,\mathbf{k}}(y) w_\alpha(y)\,dy\right)^{-1} \int_S f(y)\,p_{n,\mathbf{k}}(y) w_\alpha(y)\,dy.
$$
These operators in the univariate case and with no weight, i.e.\ $w_\alpha=1$, were introduced independently by Durrmeyer \cite{Du} and Lupa\c s \cite{Lu}; their multivariate generalization was given by Derriennic \cite{De}; and their univariate weighted form was considered by Berens and Xu \cite{Be-Xu:BD1,Be-Xu:BD2}. These operators were extensively studied by many authors and it is very difficult to summarize all the results. That is why we shall restrict our attention only to those which are directly and most closely related to the subject of the present paper. In the next section we shall recall several of their basic properties. They were proved by Ditzian \cite{Di:MBD} in the general case, and earlier by Derriennic \cite{De} and Berens and Xu \cite{Be-Xu:BD1,Be-Xu:BD2} respectively in the multivariate unweighted case and the univariate weighted case.

Ditzian \cite{Di:MBD} introduced the $K$-functional
$$
K_\alpha(f,t)_p:=\inf_{g\in C^2(S)}\left\{\|f-g\|_{p,w_\alpha} + t\,\|P_\alpha(D)g\|_{p,w_\alpha}\right\}
$$
in order to characterize the rate of approximation of the Bernstein-Durrmeyer operator in $L_{p,w_\alpha}(S)$. Here $P_\alpha(D)$ is the differential operator that is naturally associated to the multivariate Bernstein-Durrmeyer operators with the weight $w_\alpha$. It is defined by
$$
P_\alpha(D):=\sum_{\xi\in E_S} w_\alpha(x)^{-1} \frac{\partial}{\partial\xi}\,\tilde d(\xi,x)w_\alpha(x)\,\frac{\partial}{\partial\xi},
$$
where $E_S$ is the set of the directions parallel to the edges of $S$ and $\tilde d(\xi,x)$ is the distance introduced by Ditzian \cite{Di:MB}
$$
\tilde d(\xi,x):=\sup_{\substack{\lambda\ge 0\\ x+\lambda\xi\in S}} d(x,x+\lambda\xi)\sup_{\substack{\lambda\ge 0\\ x-\lambda\xi\in S}} d(x,x-\lambda\xi),
$$
as $d(x,y)$ is the Euclidean distance. 

Ditzian  \cite{Di:MBD} proved that there exist positive constants $c_1$ and $c_2$ such that for all $f\in L_p(w_\alpha)(S)$ and all $n\in\N$ there holds
\begin{equation}\label{eqchar}
c_1 K_\alpha(f,n^{-1})_p\le \|M_{n,\alpha}f-f\|_{p,w_\alpha}\le c_2 K_\alpha(f,n^{-1})_p.
\end{equation}

The direct estimate, i.e. the right-hand side inequality, was established with $c_2=2$ independently by Chen and Ditzian \cite{Ch-Di} (see also \cite[p.\ 38]{Ch-Di-Iv}) and by Berens, Schmid and Xu \cite[Theorem 2]{Be-Sc-Xu} in the unweighted case, and by Berens and Xu \cite[Theorem 3]{Be-Xu:BD1} in the univariate weighted case. A closer look at the proof of \cite[Theorem 3.3]{Di:MBD} shows that we can take $c_2$ independent of the dimension $d$ and the weight $w_\alpha$. Actually, a slight modification of this argument shows that the direct estimate holds with $c_2=2$ in the general case. More precisely, we have
\begin{equation}\label{ineqdir}
\|M_{n,\alpha}f-f\|_{p,w_\alpha}\le 2\,K_\alpha(f,n^{-1})_p.
\end{equation}
For the sake of completeness we give its proof in Section 3.

As for the converse estimate, that is, the left inequality in \eqref{eqchar}, Chen, Ditzian and Ivanov \cite[Theorems 6.1 and 6.3]{Ch-Di-Iv} established it in the unweighted case for all $d$ if $1<p<\infty$ and for $d\le 3$ if $p=1,\infty$ (a little bit weaker result was verified in the larger dimensions). Then Knoop and Zhou \cite[Theorem 3.1]{Kn-Zh} proved it for all $d$ and $1\le p\le\infty$ in the unweighted case. Both proofs give constants $c_1$ that decrease to 0 when $d$ increases. Heilmann and M. Wagner \cite[Theorem 1]{He-Wa} improved $c_1$ for $d\le 3$. Ditzian's proof of the general weighted case also yields a constant $c_1$ that decreases to 0  when $d$ or $\max_i |\alpha_i|$ increase. All these treatments are based on the quite general and efficient method developed by Ditzian and Ivanov \cite{Di-Iv}. It enables us to derive converse inequalities like the one on the left-hand side of \eqref{eqchar} by means of Voronovskaya and Bernstein-type inequalities. These inequalities are important in themselves but their consecutive application leads to decreasing $c_1$. 

The main purpose of this paper is to demonstrate that by means of the multiplier theory we can derive strong converse inequalities with better absolute constants than the methods previously used. Moreover, the arguments are very short. The first result we state contains a strong converse inequality of a form that is a combination of types B and C (according to the terminology introduced in \cite{Di-Iv}). Quite similar results were previously established by Berens and Xu \cite[Theorem 3]{Be-Xu:BD1} (see also \cite[Theorem 2]{Be-Sc-Xu}).

Set $\rho:=d+\sum_{i=1}^{d+1} \alpha_i$.

\begin{thm}\label{thm1}
Let $d\in\N$, $1\le p\le\infty$ and $w_{\alpha}$ be given by \eqref{defw} with $\alpha_i>-1$, $i=1,\dotsc,d+1$. 
Then for all $f\in L_p(w_\alpha)(S)$ and all $n\in\N$ there hold
\begin{multline*}
	K_\alpha(f,n^{-1})_p\le \left(4+\frac{2\rho}{n}\right)\big(\|M_{n,\alpha}f-f\|_{p,w_\alpha}+\|M_{2n,\alpha}f-f\|_{p,w_\alpha}\big)\\
	 +\frac{4}{n}\sum_{k=n+1}^{2n} \|M_{k,\alpha}f-f\|_{p,w_\alpha}.
\end{multline*}
\end{thm}

\begin{rem}
Let us explicitly note that the constant on the right-hand side above is asymptotically independent of any parameters unlike the strong converse inequalities obtained in \cite{Ch-Di-Iv,Di:MBD,Kn-Zh}. More precisely, if $n\ge|\rho|$, then 
\begin{multline*}
	K_\alpha(f,n^{-1})_p\le 6\big(\|M_{n,\alpha}f-f\|_{p,w_\alpha}+\|M_{2n,\alpha}f-f\|_{p,w_\alpha}\big)\\
	 +\frac{4}{n}\sum_{k=n+1}^{2n} \|M_{k,\alpha}f-f\|_{p,w_\alpha}.
\end{multline*}
However, the inequalities established in \cite{Ch-Di-Iv,Di:MBD,Kn-Zh} are of a stronger type than the one above
\end{rem}

Let us mention that the $K$-functional $K_\alpha(f,t)_p$ was characterized by a simpler one in \cite{Da-Hu-Wa} for $1<p<\infty$ (see also the references cited there). 

It seems quite plausible that the strong converse inequality in \eqref{eqchar} also holds with $c_1$, which is independent of $p$, $d$ and $w_\alpha$. We were not able to show that. However, a short multiplier argument yields a strong converse inequality of that type in a special case. It is based on a result due to H.~Pollard. Let $d=1$ and $w_\alpha=1$. Let $S_n f$ be the $n$-th partial sum of the Fourier-Legendre series of $f$. Pollard \cite{Po} proved that if $4/3<p<4$, then the operators $S_n:L_p[0,1]\to L_p[0,1]$ are uniformly bounded on $n$, that is, there exists a constant $\varsigma\ge 1$ such that
\[
\|S_n f\|_p\le\varsigma\|f\|_p,\quad f\in L_p[0,1],\ n\in\N.
\]
Here $\|\circ\|_p$ denotes the standard $L_p$-norm on the interval $[0,1]$. We will omit the subscript $\alpha$ in the notation of the $K$-functional and the Bernstein-Durrmeyer operator when $w_\alpha=1$. 

We will establish the following result.

\begin{prop}\label{pr}
Let $4/3<p<4$. Then for all $f\in L_p[0,1]$ and all $n\in\N$ there holds
\[
K(f,n^{-1})_p\le (1+2\varsigma)\,\|M_n f-f\|_p.
\]
\end{prop}
 
The contents of the paper are organized as follows. In the next section we collect the basic properties of Bernstein-Durrmeyer operator that we will use. Section 3 contains the proofs of the theorems and the proposition stated above. In the last section we discuss how the same multiplier method can be applied in the general case of weights $w_\alpha$ with $\alpha_i\ge -1/2$ for all $i$. This proof is not shorter than the ones previously used; but it has the advantage of using elementary calculus and being invariant in its technical part on the dimension---it depends only on that how large $\rho$ is.

\section{Basic properties of the Bernstein-Durrmeyer operator}

Here we shall recall the properties of the Jacobi-weighted Bernstein-Durrmeyer operator that we need (see \cite{Di:MBD}).

First of all, it is a contraction on the space $L_{p,w_\alpha}(S)$, that is,
\begin{equation}\label{bound}
\|M_{n,\alpha} f\|_{p,w_\alpha}\le \|f\|_{p,w_\alpha}.
\end{equation}

$M_{n,\alpha}$ is a self-adjoint linear operator w.r.t.\ the inner product
$$
\langle f,g \rangle_{w_\alpha}:=\int_S f(x)g(x)w_\alpha(x)\,dx.
$$
Its eigenvalues are
\begin{equation}\label{mu}
\mu_{n,\ell}:=\frac{n!}{(n-\ell)!}\,\frac{\Gamma(n+\rho+1)}{\Gamma(n+\ell+\rho+1)},\quad \ell=0,\dotsc,n,
\end{equation}
where $\Gamma$ denotes the gamma function and, to recall, we have set $\rho:=d+\sum_{i=1}^{d+1} \alpha_i$. For each $\ell$, to $\mu_{n,\ell}$ corresponds the same eigenspace for all $n$. We denote it by $V_\ell$. For $\ell\ge 1$ the space $V_\ell$ consists of those algebraic polynomials of $x_1,\dotsc,x_d$ and total degree $\ell$ that are orthogonal w.r.t.\ the above inner product to the polynomials of degree $\ell-1$. The eigenspace $V_0$, corresponding to $\mu_{n,0}=1$, consists of all constants. Now, if we denote the projections on $V_\ell$ by $\mathcal{P}_\ell$, then $M_{n,\alpha}$ can be represented in the form
\begin{equation}\label{eqexpan}
M_{n,\alpha} = \sum_{\ell=0}^n \mu_{n,\ell} \mathcal{P}_\ell.
\end{equation}

The operator $P_\alpha(D)$ is also self-adjoint and its eigenspaces coincide with those of $M_{n,\alpha}$. More precisely, there holds
\begin{equation}\label{eqP}
P_\alpha(D)P = -\ell(\ell+\rho)P,\quad P\in V_\ell,\ \ell\in\N_0.
\end{equation}

Finally, let us recall that $M_{n,\alpha}$ and $P_\alpha(D)$ commute on $C^2(S)$:
\begin{equation}\label{eqcomm}
M_{n,\alpha}P_\alpha(D) f = P_\alpha(D)M_{n,\alpha} f,\quad f\in C^2(S).
\end{equation}

\section{Proofs of the main results}

First, we will prove the direct estimate stated in \eqref{ineqdir} for the sake of completeness of the exposition.

\begin{proof}[Proof of \eqref{ineqdir}]
Z.~Ditzian's proof of the direct estimate in \eqref{eqchar}, is based on the elegant formula (see \cite[(3.3)]{Di:MBD})
\begin{equation}\label{eq4}
M_{n,\alpha}f - f =\sum_{\ell=n+1}^\infty \frac{1}{\ell(\ell+\rho)}\,P_\alpha(D)M_{\ell,\alpha}f,
\end{equation}
valid for all $f\in L_{p,w_\alpha}(S)$. Using that $M_{n,\alpha}$ is a contraction (see \eqref{bound}), we get
\begin{equation}\label{eq5}
\|M_{n,\alpha}f-f\|_{p,w_\alpha} \le 2\,\|f-g\|_{p,w_\alpha} + \|M_{n,\alpha}g-g\|_{p,w_\alpha}
\end{equation}
for any $g\in C^2(S)$. Next, we apply \eqref{bound}, \eqref{eqcomm} and \eqref{eq4} to estimate the second term on the right. Thus we get
\begin{equation}\label{eq6}
\|M_{n,\alpha}g-g\|_{p,w_\alpha} \le \sum_{\ell=n+1}^\infty \frac{1}{\ell(\ell+\rho)}\,\|P_\alpha(D) g\|_{p,w_\alpha}.
\end{equation}
It is quite straightforward, to see that
$$
\sum_{\ell=n+1}^\infty \frac{1}{\ell(\ell+\rho)} \le \frac{1}{n}.
$$
Now, substituting \eqref{eq6} in \eqref{eq5} and taking an infimum on $g\in C^2(S)$, we arrive at
\begin{equation*}
\|M_{n,\alpha}f-f\|_{p,w_\alpha}\le 2\,K_\alpha(f,n^{-1})_p.
\end{equation*}
Thus the first inequality in \eqref{ineqdir} is verified; the second one is trivial.
\end{proof}

\begin{proof}[Proof of \thmref{thm1}]
The proof is a modification of a very short argument due to Berens and Xu (see \cite[Theorem 3]{Be-Xu:BD1}). Set
\[
g_n:=\frac{1}{t_n}\sum_{k=n+1}^{2n} \frac{M_{k,\alpha}f}{k(k+\rho)},\quad t_n:=\sum_{k=n+1}^{2n} \frac{1}{k(k+\rho)}.
\]
Clearly, $g_n\in C^2(S)$ for all $n\in\N$ and then
\begin{equation}\label{eq20}
K_\alpha(f,n^{-1})_p\le \|f-g_n\|_{p,w_\alpha} + \frac{1}{n}\,\|P_\alpha(D) g_n\|_{p,w_\alpha}.
\end{equation}
We estimate the first term on the right above by means of
\begin{equation}\label{eq21}
\begin{split}
	\|f-g_n\|_{p,w_\alpha}&=\left\|f-\frac{1}{t_n}\sum_{k=n+1}^{2n} \frac{M_{k,\alpha}f}{k(k+\rho)}\right\|_{p,w_\alpha}\\
	&\le \frac{1}{t_n}\sum_{k=n+1}^{2n} \frac{\|M_{k,\alpha}f-f\|_{p,w_\alpha}}{k(k+\rho)}\\
	&\le \frac{4}{n}\sum_{k=n+1}^{2n} \|M_{k,\alpha}f-f\|_{p,w_\alpha}.
	\end{split}
\end{equation}
In order to estimate the second term on the right in \eqref{eq20}, we apply \eqref{eqexpan} and \eqref{eqP} to get the representation
\begin{equation*}
	P_\alpha(D) g_n	= -\frac{1}{t_n}\,\sum_{k=n+1}^{2n}\sum_{\ell=0}^k \frac{\ell(\ell+\rho)}{k(k+\rho)}\,\mu_{k,\ell} \mathcal{P}_\ell.
\end{equation*}
Next, we take into account the remarkable property of the multipliers $\mu_{n,\ell}$
\[
\mu_{k,\ell}-\mu_{k-1,\ell}=\frac{\ell(\ell+\rho)}{k(k+\rho)}\,\mu_{k,\ell}
\]
to arrive at the formula
\begin{align*}
	P_\alpha(D) g_n &= \frac{1}{t_n}\,\sum_{k=n+1}^{2n}(M_{k-1,\alpha}f-M_{k,\alpha}f)\\
	&= \frac{1}{t_n}\,(M_{n,\alpha}f-M_{2n,\alpha}f).
\end{align*}
Consequently,
\begin{multline}\label{eq22}
\frac{1}{n}\,\|P_\alpha(D) g_n\|_{p,w_\alpha}\\
\le \left(4+\frac{2\rho}{n}\right)\big(\|M_{n,\alpha}f-f\|_{p,w_\alpha} + \|M_{2n,\alpha}f-f\|_{p,w_\alpha}\big).
\end{multline}
Combining \eqref{eq20}-\eqref{eq22}, we complete the proof of the theorem.
\end{proof}

Let us proceed to the proof of the converse inequality in \propref{pr}. The method we use is quite straightforward. It is based entirely on standard techniques in the multiplier theory and orthogonal series expansions. We will present it in the general case of the multivariate Bernstein-Durrmeyer operator on the simplex. The method is based on constructing a family of uniformly bounded operators $\mathcal{Q}_n$ such that 
$$
\frac{1}{n}\,P_\alpha(D)M_{n,\alpha}^m f=\mathcal{Q}_n(M_{n,\alpha}f-f)
$$
with some fixed $m\in\N$. Then the strong one-term converse inequality in \eqref{eqchar} easily follows from 
\begin{align*}
		&K_\alpha(f,n^{-1})_p \le \|M_{n,\alpha}^m f-f\|_{p,w_\alpha} + \frac{1}{n}\,\|P_\alpha(D)M_{n,\alpha}^m f\|_{p,w_\alpha}\\
		&\qquad =\|(M_{n,\alpha}^{m-1}+M_{n,\alpha}^{m-2}+\dotsb+I)(M_{n,\alpha} f-f)\|_{p,w_\alpha} + \|\mathcal{Q}_n(M_{n,\alpha}f-f)\|_{p,w_\alpha}\\
		&\qquad\le (m+q)\,\|M_{n,\alpha} f-f\|_{p,w_\alpha},
\end{align*}
where $I$ denotes the identity and $q>0$ is such that $\|\mathcal{Q}_n F\|_{p,w_\alpha}\le q\|F\|_{p,w_\alpha}$ for all $F\in L_{p,w_\alpha}(S)$ and $n\in\N$.

That approach to converse inequalities has been applied before (see e.g.\ \cite[(2.13)]{Di-Iv}, and also cf.\ \cite[p.\ 32]{Be-Xu:BD1}). The proof of the direct inequality, we recalled above, was realized in a similar way (see \eqref{eq4}). There is a general comparison principle that underlies this technique. It was formulated independently, in two different settings, by Shapiro \cite{Sh:CP} (see also \cite[Section 9.4]{Sh:AT}) and Trigub \cite[\S\,4]{Tr:Sum} and \cite[\S\,4]{Tr:Abs} (see also \cite[Chapter 7]{Tr-Be} and \cite[p.\ 4]{Tr:K}. The author tried to present systematically that method of verifying direct and converse estimates in terms of $K$-functionals in \cite{Dr} (see also the references cited there).

The earlier proofs of the converse inequality of the type given in \eqref{eqchar} for the Bernstein-Durrmeyer operator also employed orthogonal expansions, but in a lesser degree and within the framework in \cite{Di-Iv}. Berens and Xu \cite{Be-Xu:BD1} also extensively used multiplier techniques (see also \cite[Theorem 2]{Be-Sc-Xu}). 

\begin{proof}[Proof of \propref{pr}]
Let us begin with several observations valid in the general multivariate weighted case. They will be useful for our discussion in the next section. 

We first note that \eqref{eqexpan} and \eqref{eqP} yield
$$
P_\alpha(D)M_{n,\alpha} f = -\sum_{\ell=1}^n \ell(\ell+\rho) \mu_{n,\ell}\mathcal{P}_\ell f.
$$
We introduce the linear operator on $L_{p,w_\alpha}(S)$
$$
Q_n f:= \sum_{\ell=1}^n \nu_{n,\ell} \mathcal{P}_\ell f,
$$
where
\begin{equation}\label{eqnu}
\nu_{n,\ell}:=\frac{\ell(\ell+\rho)\,\mu_{n,\ell}}{n(1-\mu_{n,\ell})}.
\end{equation}
Note that $\mu_{n,\ell}<1$ for $\ell=1,2,\dotsc,n$. With that operator we have
$$
\frac{1}{n}\,P_\alpha(D)M_{n,\alpha} f= Q_n(M_{n,\alpha} f-f).
$$
Thus to establish a one-term strong converse inequality, it is enough to show that
$$
\|Q_n f\|_{p,w_\alpha}\le c\,\|f\|_{p,w_\alpha}
$$
for all $f\in L_{p,w_\alpha}(S)$ and $n\in\N$.

After this general remark, we proceed to the proof of the proposition. Now, $S_n f$ coincide with the $n$th partial sum of the orthogonal expansion of $f$ on $\mathcal{P_\ell}$, that is,
\[
S_n f:=\sum_{\ell=0}^n P_\ell f.
\] 

We use the representation
\begin{equation*}
Q_n f =\sum_{\ell=1}^{n-1} (\nu_{n,\ell}-\nu_{n,\ell+1})S_\ell f + \nu_{n,n}S_n f - \nu_{n,1}S_0 f.
\end{equation*}
In \lemref{l1} below we will show that $\nu_{n,\ell}-\nu_{n,\ell+1}>0$ for all $\ell$. Then, taking also into account that the $\nu$'s are positive and $\nu_{n,1}=1$, we deduce the estimate
\begin{align*}
	\|Q_n f\|_{p,w_\alpha}&\le \varsigma\left(\sum_{\ell=1}^{n-1} (\nu_{n,\ell}-\nu_{n,\ell+1}) + \nu_{n,n} + \nu_{n,1}\right)\|f\|_{p,w_\alpha}\\
	&\le 2\varsigma\nu_{n,1}\|f\|_{p,w_\alpha}=2\varsigma\,\|f\|_{p,w_\alpha};
\end{align*}
hence the assertion of the proposition follows.
\end{proof}

\begin{lem}\label{l1}
Let $\rho>-1$. For $\nu_{n,\ell}$ defined in \eqref{eqnu} there holds
\begin{equation}\label{eq2}
\nu_{n,\ell}>\nu_{n,\ell+1},\quad \ell=1,\dotsc,n-1,\ n=2,3,\dotsc.
\end{equation}
\end{lem}

\begin{proof}
Relation \eqref{eq2} is equivalent to
\[
\frac{1-\mu_{n,\ell}}{\ell(\ell+\rho)\,\mu_{n,\ell}}<\frac{1-\mu_{n,\ell+1}}{(\ell+1)(\ell+\rho+1)\,\mu_{n,\ell+1}},
\]
which can be written in the form
\[
\frac{1}{\ell(\ell+\rho)\,\mu_{n,\ell}}-\frac{1}{\ell(\ell+\rho)}<\frac{1}{(\ell+1)(\ell+\rho+1)\mu_{n,\ell+1}}-\frac{1}{(\ell+1)(\ell+\rho+1)}.
\]
We group the terms with $\mu$'s on the left-hand side and those without on the right-hand side, and substitute the value of the $\mu$'s given in \eqref{mu}. After straightforward calculations, using that $\rho>-1$ and
\begin{equation}\label{eq9}
\Gamma(n+\ell+\rho+2)=(n+\ell+\rho+1)\Gamma(n+\ell+\rho+1),
\end{equation}
which follows from $\Gamma(z+1)=z\,\Gamma(z)$, $z>0$, we deduce that \eqref{eq2} is equivalent to
\[
(n-\ell-1)!\,\Gamma(n+\ell+\rho+1)[n-\ell(\ell+\rho+1)]<n!\,\Gamma(n+\rho+1)
\]
for $\ell=1,\dotsc,n-1,\ n=2,3,\dotsc$. To verify this inequality, we shall show that the quantity on the left-hand side is decreasing on $\ell$ and it is valid for $\ell=1$. The latter is a matter of a direct check---it reduces to $(\rho+1)(\rho+2)>0$. To verify the former, we set
\[
\xi_{n,\ell}:=(n-\ell-1)!\,\Gamma(n+\ell+\rho+1)[n-\ell(\ell+\rho+1)].
\]
To see that
\begin{equation}\label{eq3}
\xi_{n,\ell}>\xi_{n,\ell+1},\quad \ell=1,2,\dotsc,n-2,\ n=3,4,\dotsc,
\end{equation}
we again apply \eqref{eq9} to deduce that \eqref{eq3} is equivalent to
\[
(n-\ell-1)[n-\ell(\ell+\rho+1)]>(n+\ell+\rho+1)[n-(\ell+1)(\ell+\rho+2)].
\]
Now, direct computations yield
\begin{multline*}
(n-\ell-1)[n-\ell(\ell+\rho+1)]-(n+\ell+\rho+1)[n-(\ell+1)(\ell+\rho+2)]\\
=(\ell+1)(\ell+\rho+1)(2\ell+\rho+2)>0,
\end{multline*}
which verifes \eqref{eq3} and completes the proof of the lemma.
\end{proof}

\section{An extension}

In this section we will demonstrate that the multiplier method can be used to verify the one-term converse inequality in \eqref{eqchar} in a more general situation than the one considered in \propref{pr}.

To this end, we represent $Q_n$ as a linear combination of the Ces\`aro means of the partial sums of the orthogonal expansion of $f$ on $\mathcal{P_\ell}$ (see \cite[Theorem 3.2]{Bu-Ne-Tr}). We set
\[
\widetilde S_n f:=\frac{1}{n+1}\sum_{k=0}^n S_k f,
\]
where
\[
S_k f:=\sum_{\ell=0}^k P_\ell f.
\] 

Then we have
\begin{multline}\label{eq24}
	Q_n f =\sum_{\ell=1}^{n-2} (\ell+1)(\nu_{n,\ell+2}-2\nu_{n,\ell+1}+\nu_{n,\ell})\widetilde S_\ell f\\
	+n(\nu_{n,n-1}-2\nu_n)\widetilde S_{n-1}f +(n+1)\nu_{n,n}\widetilde S_n f +(\nu_{n,2}-2\nu_{n,1})\widetilde S_0 f.
\end{multline}
As usually, if the range of summation is empty, we set the sum to be equal to zero.

Dai and Xu \cite[Theorem 2.8 with $\delta=1$]{Da-Xu:BC2} (or see \cite[Theorem 13.4.4]{Da-Xu:AT}, as we also apply the Riesz-Thorin interpolation theorem) showed that if $1\le p\le\infty$, $\alpha_i\ge-1/2$, $i=0,\dotsc,d+1$, and $\rho - \min_i\alpha_i<3/2$, then the Ces\`aro means are uniformly bounded on $n$, i.e.~there exists a constant $\kappa$ such that
\begin{equation}\label{eq39}
\|\widetilde S_n f\|_{p,w_\alpha}\le\kappa\,\|f\|_{p,w_\alpha},\quad f\in L_p(w_\alpha)(S),\ n\in\N.
\end{equation}

\lemref{l1} yeilds $\nu_{n,2}\le\nu_{n,1}=1$. Then we have by \eqref{eq24} and \eqref{eq39}
\begin{multline*}
	\|Q_n f\|_{p,w_\alpha} \le\kappa\Bigg(\sum_{\ell=1}^{n-2} (\ell+1)|\nu_{n,\ell+2}-2\nu_{n,\ell+1}+\nu_{n,\ell}|\\
	+(4n+1)\nu_{n,n-1} +3\Bigg)\|f\|_{p,w_\alpha}.
\end{multline*}
We will prove that
\[
\sum_{\ell=1}^{n-2} (\ell+1)|\nu_{n,\ell+2}-2\nu_{n,\ell+1}+\nu_{n,\ell}|\le c
\]
and
\[
n\,\nu_{n,n-1}\le c.
\]
Above and henceforward, $c$ denotes a positive constant, not necessarily the same at each occurrence, whose value is independent of $n$.

Thus we will have shown that if $1\le p\le\infty$, $\alpha_i\ge -1/2$, $i=1,\dotsc,d+1$, and
\[
d+\sum_{i=1}^{d+1} \alpha_i - \min_{1\le i\le d+1}\alpha_i<\frac{3}{2},
\]
then for all $f\in L_p(w_\alpha)(S)$ and all $n\in\N$ there holds
\[
K_\alpha(f,n^{-1})_p\le c\,\|M_{n,\alpha}f-f\|_{p,w_\alpha}.
\]
In order to treat the general case, we can still apply the same method but use Ces\`aro means of higher order (see \cite[Theorem 7.1]{Bu-Ne-Tr:2} or \cite[Theorem 3.3]{Tr:M}). Their uniform boundedness was established by Dai and Xu \cite{Da-Xu:BC2} (or see \cite[Theorems 13.2.7 and 13.4.6]{Da-Xu:AT}).

We proceed to establishing the auxiliary results.

We set for $\tau\in (0,n]$
\[
\mu_n(\tau):=\frac{\Gamma(n+1)\Gamma(n+\rho+1)}{\Gamma(n-\tau+1)\Gamma(n+\tau+\rho+1)},\quad
\nu_n(\tau):=\frac{\tau(\tau+\rho)\,\mu_n(\tau)}{n(1-\mu_n(\tau))}.
\]
We will make use of the following formula of the derivative of the gamma function
\[
\Gamma'(z)=\Gamma(z)\psi(z),
\]
where $\psi(z)$ is the digamma function, defined as the logarithmic derivative of the gamma function
\[
\psi(z):=\frac{\Gamma'(z)}{\Gamma(z)}.
\]
We have 
\begin{equation}\label{eqrepr}
\mu_n'(\tau)=-\mu_n(\tau)C_n(\tau),
\end{equation}
where
\[
C_n(\tau):=\psi(n+\tau+\rho+1)-\psi(n-\tau+1).
\]

We will use the following estimates.

\begin{lem}\label{l5}
Let $\rho\ge 0$. Then:
\begin{align}
	&C_n(\tau)\le \frac{2\tau+\rho}{n-\tau},\quad \tau\in (0,n);\label{l5a}\\
	&C_n(\tau)\ge\frac{2\tau+\rho}{2(n-\tau+1)},\quad \tau\in (0,(n-\rho)/3),\ n>\rho;\label{l5b}\\
	&C'_n(\tau)\le \frac{2n+\rho}{(n+\tau+\rho)(n-\tau)},\quad \tau\in (0,n);\label{l5c}\\
	&C'_n(\tau)\ge \frac{2n+\rho+2}{(n+\tau+\rho+1)(n-\tau+1)},\quad \tau\in (0,n);\label{l5d}\\
	&C''_n(\tau)\ge \frac{2(2\tau+\rho-1)(2n+\rho+1)}{(n+\tau+\rho)^2(n-\tau+1)^2},\quad \tau\in (0,n).\label{l5e}
\end{align}
\end{lem}

\begin{proof}
As is known,
\begin{equation}\label{eq32}
\psi(x)=-\gamma-\frac{1}{x}+\sum_{k=1}^\infty \frac{x}{k(k+x)},\quad x>0,
\end{equation}
where $\gamma$ is Euler's constant. Therefore
\begin{equation}\label{eq33}
C_n(\tau)=(2\tau+\rho)\sum_{k=1}^\infty \frac{1}{(n-\tau+k)(n+\tau+\rho+k)}.
\end{equation}
Interpreting the sum above as a Darboux sum, we arrive at the estimates
\begin{equation}\label{eq34}
\log\left(1+\frac{2\tau+\rho}{n-\tau+1}\right)\le C_n(\tau)\le \log\left(1+\frac{2\tau+\rho}{n-\tau}\right).
\end{equation}
To complete the proof of the first two estimates, it remains to take into account the inequalities
\begin{align*}
	&\log(1+x)\le x,\quad x\in\R,\\
	&\log(1+x)\ge x-\frac{x^2}{2}\ge\frac{x}{2},\quad x\in [0,1].
\end{align*}

In order to estimate the derivatives of $C_n$, we use that for $m\ge 1$ we have
\[
\psi^{(m)}(x)=(-1)^{m+1} m!\sum_{k=0}^\infty \frac{1}{(x+k)^{m+1}},\quad x>0.
\]
Therefore
\begin{gather}
	\frac{1}{x}\le\psi'(x)\le\frac{1}{x-1};\label{eq35}\\
	-\frac{2}{(x-1)^2}\le\psi''(x)\le -\frac{2}{x^2};\label{eq35a}\\
\end{gather}
for $x>1$. These inequalities directly yield \eqref{l5c}-\eqref{l5e}.
\end{proof}

\begin{lem}\label{l6}
Let $\rho\ge 0$, $b>0$ and $0<\delta\le 1$. Let also $n\in\N$ be such that $n\ge 3$ and $1\le\sqrt{bn}\le n-1$. Then
\begin{align}
	&n^2\nu_{n,\ell}\le c,\quad \delta n\le\ell\le n,\label{l60}\\
	&\tau|\nu_n'(\tau)|\le c,\quad\tau\in [1,n-1],\label{l6a}
	\intertext{and}
	&\tau^2|\nu_n''(\tau)|\le c,\quad\tau\in [1,\sqrt{bn}],\label{l6b}
\end{align}
where the constant $c$ is independent of $n$.
\end{lem}

\begin{proof}
First, we estimate from below the difference $1-\mu_{n,\ell}$.

By means of the property $\Gamma(z+1)=z\Gamma(z)$, $z>0$, we represent $\mu_{n,\ell}$ in the form
\[
\mu_{n,\ell}=\frac{n(n-1)\dotsb(n-\ell+1)}{(n+\rho+1)(n+\rho+2)\dotsb(n+\rho+\ell)}.
\]
Consequently,
\[
1-\mu_{n,\ell}\ge \frac{(n+\rho+1)(n+\rho+2)\dotsb(n+\rho+\ell)-n^\ell}{(n+\rho+1)(n+\rho+2)\dotsb(n+\rho+\ell)}.
\]
We expand the numerator, take into account that $\rho\ge 0$, and use the well-known formulas for sums of powers of consecutive positive integers, to arrive at the estimate
\[
(n+\rho+1)(n+\rho+2)\dotsb(n+\rho+\ell)-n^\ell\ge c(\ell^2 n^{\ell-1}+\ell^6 n^{\ell-3}).
\]
Hence we get the inequalities
\begin{align}
	&1-\mu_{n,\ell}\ge \frac{c\,\ell^2 n^{\ell-1}}{(n+\rho+1)(n+\rho+2)\dotsb(n+\rho+\ell)}\label{eq31}
	\intertext{and}
	&1-\mu_{n,\ell}\ge \frac{c\,\ell^6 n^{\ell-3}}{(n+\rho+1)(n+\rho+2)\dotsb(n+\rho+\ell)}\label{eq31a}
\end{align}
for $3\le\ell\le n$.

Inequality \eqref{l60} for $\ell\ge 3$ follows directly from \eqref{eq31a} and $\ell\ge\delta n$:
\[
	n^2\nu_{n,\ell}\le c\,\frac{n^{\ell+2}}{\ell^4 n^{\ell-2}}\le c.
\]
For $\ell=1,2$ \eqref{l60} is trivial.

We proceed to the second assertion of the lemma. Making use of \eqref{eqrepr}, we arrive at
\begin{equation}\label{eq34a}
\tau\nu_n'(\tau)=\frac{\tau(2\tau+\rho)\,\mu_n(\tau)}{n(1-\mu_n(\tau))}-\frac{\tau^2(\tau+\rho)\,\mu_n(\tau)\,C_n(\tau)}{n(1-\mu_n(\tau))^2}.
\end{equation}

The function $\mu_n(\tau)$ is monotone decreasing on $\tau$ for each fixed $n$. For the rest of the proof let $\ell\in\{1,\dotsc,n-2\}$ be such that $\ell\le\tau\le\ell+1$. Then
\begin{align}
	&\mu_n(\tau)\le\mu_{n,\ell},\label{eq29}\\
	&1-\mu_n(\tau)\ge 1-\mu_{n,\ell}\label{eq30}.
\end{align}
These two inequalities, the property $\Gamma(z+1)=z\Gamma(z)$, $z>0$, and \eqref{eq31} imply the following estimate of the first term on the right in \eqref{eq34a} 
\begin{equation}\label{eq36}
\begin{split}
0\le\frac{\tau(2\tau+\rho)\,\mu_n(\tau)}{n(1-\mu_n(\tau))}
&\le \frac{(\ell+1)(2\ell+\rho+2)\,\mu_{n,\ell}}{n(1-\mu_{n,\ell})}\\
&\le c\,\frac{(\ell+1)(2\ell+\rho+2)}{\ell^2}\,\frac{n!}{n^\ell (n-\ell)!}\\
&\le c,\quad \tau\in [1,n-1].
\end{split}
\end{equation}

To estimate the second term we argue in a similar way, as we also use \eqref{l5a}. We have
\begin{equation*}
\begin{split}
\frac{\tau^2(\tau+\rho)\,\mu_n(\tau)\,nC_n(\tau)}{(n(1-\mu_n(\tau)))^2}
&\le c\,\frac{(\ell+1)^2 (\ell+\rho+1)(2\ell+\rho+2)}{\ell^4}\\
&\qquad \times\frac{n!\,(n+\rho+1)\dotsb(n+\rho+\ell)}{n^{2\ell-1}(n-\ell)!\,(n-\ell-1)}\\
&\le c\left(1-\frac{1}{n}\right)\dotsb\left(1-\frac{\ell-2}{n}\right)\!\left(1+\frac{\rho+1}{n}\right)\dotsb\left(1+\frac{\rho+\ell}{n}\right)\\
&\le c\prod_{i=1}^{\ell-2}\left(1-\frac{i}{n}\right)\left(1+\frac{i+\rho}{n}\right).
\end{split}
\end{equation*}
As usually, we set an empty product to be equal to $1$.

Next, we take into account that
\begin{equation}\label{eq54}
	\left(1-\frac{i}{n}\right)\left(1+\frac{i+\rho}{n}\right)
	=1-\frac{i^2}{n^2}+\frac{\rho}{n}\left(1-\frac{i}{n}\right)\le 1+\frac{\rho}{n}
\end{equation}
and the inequality $(1+\rho/n)^n\le e^\rho$ to deduce
\begin{equation}\label{eq37}
0\le\frac{\tau^2(\tau+\rho)\,\mu_n(\tau)\,C_n(\tau)}{n(1-\mu_n(\tau))^2}\le c,\quad \tau\in [1,n-1].
\end{equation}

Relations \eqref{eq34a}, \eqref{eq36} and \eqref{eq37} imply the second inequality in the lemma.

In order two prove the last assertion of the lemma, we use the representation
\begin{multline}\label{eq38}
\nu''(\tau)=\frac{2\mu_n(\tau)}{n(1-\mu_n(\tau))}-\frac{2(2\tau+\rho)\,\mu_n(\tau)\,C_n(\tau)}{n(1-\mu_n(\tau))^2}\\
-\frac{\tau(\tau+\rho)\,\mu_n(\tau)\,C_n'(\tau)}{n(1-\mu_n(\tau))^2}
+\frac{\tau(\tau+\rho)(1+\mu_n(\tau))\,\mu_n(\tau)\,C_n(\tau)^2}{n(1-\mu_n(\tau))^3}.
\end{multline}

Just similarly to \eqref{eq36} and \eqref{eq37}, we establish
\begin{gather}
0\le\frac{\tau^2\mu_n(\tau)}{n(1-\mu_n(\tau))}\le c,\label{eq50}\\
0\le\frac{\tau^2(2\tau+\rho)\mu_n(\tau)nC_n(\tau)}{(n(1-\mu_n(\tau)))^2}\le c\label{eq51}
\end{gather}
for $\tau\in [1,n-1]$.

Again, similarly to the proof of \eqref{eq37}, but this time using \eqref{l5c}, we get
\begin{align*}
\frac{\tau^3(\tau+\rho)\,\mu_n(\tau)\,C_n'(\tau)}{n(1-\mu_n(\tau))^2}
&\le c\,\frac{(\ell+1)^3 (\ell+\rho+1)}{\ell^4}\\
&\qquad \times\frac{(2n+\rho)\,n!\,(n+\rho+1)\dotsb(n+\rho+\ell-1)}{n^{2\ell-1}(n-\ell)!\,(n-\ell-1)}\\
&\le c\prod_{i=1}^{\ell-2}\left(1-\frac{i}{n}\right)\left(1+\frac{i+\rho}{n}\right)\le c.
\end{align*}
Consequently,
\begin{equation}\label{eq52}
0\le \frac{\tau^3(\tau+\rho)\,\mu_n(\tau)\,C_n'(\tau)}{n(1-\mu_n(\tau))^2}\le c,\quad \tau\in [1,n-1].
\end{equation}

In order to estimate the last term in the representation of $\nu''_n$ we use \eqref{l5a} and $\mu_{n,\ell}\le 1$ to deduce
\begin{align*}
\frac{\tau^3(\tau+\rho)(1+\mu_n(\tau))\,\mu_n(\tau)\,C_n(\tau)^2}{n(1-\mu_n(\tau))^3}
&\le c\,\frac{(\ell+1)^3 (\ell+\rho+1)(2\ell+\rho+2)^2}{\ell^6}\\
&\qquad \times\frac{n!\,(n+\rho+1)^2\dotsb(n+\rho+\ell)^2}{n^{3\ell-2}(n-\ell)!\,(n-\ell-1)^2}\\
&\le c\prod_{i=1}^{\ell-3}\left(1-\frac{i}{n}\right)\left(1+\frac{i+\rho}{n}\right)^2.
\end{align*}
It remains to observe that, by virtue of \eqref{eq54} and the inequality $(1+\rho/n)^n\le e^\rho$, we have
\begin{align*}
	\prod_{i=1}^{\ell-3}\left(1-\frac{i}{n}\right)\left(1+\frac{i+\rho}{n}\right)^2
	&\le c \left[\left(1+\frac{\ell+\rho}{n}\right)^n\right]^{\ell/n}\le c\,e^{\ell^2/n}\le c.
\end{align*}
\end{proof}

\begin{lem}\label{l3}
Let $\rho\ge 0$. There holds
\[
\ell(\nu_{n,\ell}-\nu_{n,\ell+1})\le c,\quad \ell=1,\dotsc,n-1,
\]
where the constant $c$ is independent of $n$.
\end{lem}

\begin{proof}
The inequality follows readily from \eqref{l60} for $\ell=n-1$. Let $\ell\le n-2$. Then, by virtue of \eqref{l6a}, we have
\begin{equation*}
\begin{split}
\ell(\nu_{n,\ell}-\nu_{n,\ell+1})&=-\ell\int_\ell^{\ell+1} \nu_n'(\tau)\,d\tau\\
&\le \sup_{1\le\tau\le n-1}\left|\tau\nu_n'(\tau)\right|\le c.
\end{split}
\end{equation*}
\end{proof}

\begin{lem}\label{l4}
Let $\rho\ge 0$. There holds
\[
\sum_{\ell=1}^{n-2} (\ell+1)|\nu_{n,\ell+2}-2\nu_{n,\ell+1}+\nu_{n,\ell}|\le c,
\]
where the constant $c$ is independent of $n$.
\end{lem}

\begin{proof}
Clearly, it is sufficient to verify the lemma for large $n$. Its assertion for $n\le n_0$, where $n_0\in\N$ is fixed, is trivial.

We split the sum into four parts:
\begin{align*}
	&1\le \ell\le\sqrt{an}-2,\\
	&\sqrt{an}-2<\ell\le\sqrt{bn},\\
	&\sqrt{bn}<\ell\le\frac{n}{4},\\
	&\frac{n}{4}<\ell\le n-2,	
\end{align*}
where $0<a<b$ will be fixed in appropriate way to be indicated in the course of the proof. We denote these parts with $\Sigma_i$, $i=1,\dotsc,4$, respectively.

As is known
\[
\nu_{n,\ell+2}-2\nu_{n,\ell+1}+\nu_{n,\ell}=\int_\ell^{\ell+2} M(\tau-\ell) \nu_n''(\tau)\,d\tau,\quad \ell=1,\dotsc,n-2,
\]
where
\[
M(\tau):=\begin{cases}
\tau, &0\le \tau\le 1,\\
2-\tau, & 1\le \tau\le 2.
\end{cases}
\]

By virtue of \eqref{l6b}, we have
\begin{align*}
	\Sigma_2&:=\sum_{\sqrt{an}-2<\ell\le\sqrt{bn}} (\ell+1)|\nu_{n,\ell+2}-2\nu_{n,\ell+1}+\nu_{n,\ell}|\\
	&\le c\int_{\sqrt{an}-2}^{\sqrt{bn}+2} \tau|\nu_n''(\tau)|\,d\tau\le c.
\end{align*}

Let $m_n$ be the integer part of $n/4$. We apply \eqref{l60} to get
\[
	\Sigma_4:=\sum_{n/4<\ell\le n-2} (\ell+1)|\nu_{n,\ell+2}-2\nu_{n,\ell+1}+\nu_{n,\ell}|\\
	\le c\,n^2\nu_{n,m_n}\le c.
\]

We proceed to estimating $\Sigma_3$. Let $\sqrt{bn}\le\tau\le n/4+2$. Let $n$ be so large that we have $n/4+2\le(n-\rho)/3$. We will show that if $b$ is fixed large enough, then $\nu_n''(\tau)>0$ for all large $n$. Hence $\nu_{n,\ell+2}-2\nu_{n,\ell+1}+\nu_{n,\ell}\ge 0$ if $\sqrt{bn}<\ell\le n/4-2$. Let $\ell_n$ be the smallest integer greater than $\sqrt{bn}$. Then, by virtue also of Lemmas \ref{l1} and \ref{l3}, we deduce that
\begin{align*}
	\Sigma_3&:=\sum_{\sqrt{bn}<\ell\le n/4} (\ell+1)|\nu_{n,\ell+2}-2\nu_{n,\ell+1}+\nu_{n,\ell}|\\
	&=\ell_n(\nu_{n,\ell_n}-\nu_{n,\ell_n+1})+\nu_{n,\ell_n}-(m_n+1)\nu_{n,m_n}+m_n\nu_{n,m_n+1}\\
	&\le c.
\end{align*}

Thus to complete the proof of the estimate of $\Sigma_3$ it remains to show $\nu_{n,\ell+2}-2\nu_{n,\ell+1}+\nu_{n,\ell}\ge 0$ if $\sqrt{bn}\le\ell\le n/4-2$ for all $n$ large enough an appropriately fixed $b$. By \eqref{eq38} we have
\[
\nu_n''(\tau)=\frac{\mu_n(\tau)}{n(1-\mu_n(\tau))^3}N(\tau),
\]
where we have set
\begin{multline*}
	N(\tau):=2(1-\mu_n(\tau))^2-2(2\tau+\rho)C_n(\tau)(1-\mu_n(\tau))\\
		-\tau(\tau+\rho)C_n'(\tau)(1-\mu_n(\tau))+\tau(\tau+\rho)C_n^2(\tau)(1+\mu_n(\tau)).
\end{multline*}
By virtue of \lemref{l5}, we arrive at the estimate
\begin{multline}
	N(\tau)\ge \mu_n(\tau)
	\left(\frac{2(2\tau+\rho)^2}{n-\tau}+\frac{\tau(\tau+\rho)(2n+\rho)}{(n+\tau+\rho)(n-\tau)}
	+\frac{\tau(\tau+\rho)(2\tau+\rho)^2}{4(n-\tau+1)^2}\right)\\
	-\frac{2(2\tau+\rho)^2}{n-\tau}-\frac{\tau(\tau+\rho)(2n+\rho)}{(n+\tau+\rho)(n-\tau)}
	+\frac{\tau(\tau+\rho)(2\tau+\rho)^2}{4(n-\tau+1)^2}.
\end{multline}
In order to show that $N(\tau)>0$ it is enough to prove that the quantity on the right-hand side of the last relation is positive. Using that $n-\tau+1<n+\tau+\rho$, we see that this follows from
\begin{multline*}
	\mu_n(\tau)[8(2\tau+\rho)^2(n+\tau+\rho)^2+4\tau(\tau+\rho)(2n+\rho)(n+\tau+\rho)+\tau(\tau+\rho)(2\tau+\rho)^2(n-\tau)]\\
	>8(2\tau+\rho)^2(n+\tau+\rho)^2+4\tau(\tau+\rho)(2n+\rho)(n+\tau+\rho)-\tau(\tau+\rho)(2\tau+\rho)^2(n-\tau).
\end{multline*}
To complete the proof it remains to observe that if $b$ is fixed large enough, then the quantity on the right-hand side of the inequality above is negative for large $n$. To see this, we observe that the sum of the terms in the polynomial on the variables $\tau$ and $n$ on the right-hand side that determine its sign for large $\tau$ and $n$ is
\[
40n^2\tau^2+72n\tau^3-8\rho n\tau^3-4n\tau^4+4\tau^5.
\]
Since
\[
40n^2\tau^2+72n\tau^3-8\rho n\tau^3-4n\tau^4+4\tau^5\le 4\tau^2(10n^2+18n\tau-n\tau^2+\tau^3),
\]
to complete the proof it is sufficient to show that
\begin{equation*}
10n^2+18n\tau-n\tau^2+\tau^3<0
\end{equation*}
if $\sqrt{bn}\le\tau\le n/4$ with an appropriately fixed $b$. But this readily becomes clear from the estimate
\begin{align*}
	10n^2+18n\tau-n\tau^2+\tau^3
	&\le 10n^2+\frac{9}{2}n^2-n\tau^2+\frac{1}{4}n\tau^2\\
	&\le \frac{29}{2}n^2-\frac{3b}{4} n^2.
\end{align*}

To estimate $\Sigma_1$ we use similar but more lengthy considerations than those for $\Sigma_3$. They are based on the inequalities stated in \lemref{l5} as we have to use instead of \eqref{l5b} its refinement that follows from $\log(1+x)\ge x-x^2/2$. This time we show that there exists $a\in (0,1)$ such that $N'(\tau)<0$ at least for large $n$ if $1\le\tau\le\sqrt{an}$; hence $N(\tau)\le N(1)<0$. Consequently, $\nu_{n,\ell+2}-2\nu_{n,\ell+1}+\nu_{n,\ell}\le 0$ if $1\le\ell\le\sqrt{an}-2$ and $n$ is large. 
\end{proof}

\bigskip

\textbf{Acknowledgments.} Supported by grant 80-10-13/2018 of the Research Fund of the University of Sofia. I am thankful to Professor Dany Leviatan for his helpful comments that improved the presentation of the results.

\bigskip
\noindent
\begin{footnotesize}
\begin{tabular}{ll}
Borislav R. Draganov&
\\
Dept. of Mathematics and Informatics&
Inst. of Mathematics and Informatics\\
University of Sofia&
Bulgarian Academy of Sciences\\
5 James Bourchier Blvd.&
bl. 8 Acad. G. Bonchev Str.\\
1164 Sofia&
1113 Sofia\\
Bulgaria&
Bulgaria\\
bdraganov@fmi.uni-sofia.bg&
\\
\end{tabular}
\noindent
\end{footnotesize}

\end{document}